
\magnification=1200
\hsize=15truecm
\vsize=24truecm
\centerline{\bf Energy of embedded surfaces invariant
under M\oe bius transformations, addendum}
\medskip
\centerline{Stefano DEMICHELIS}
\medskip
{\bf Abstract:}In a previous preprint we defined an energy associated to
 every embedding
of a surface in $R^n$ or $S^n$. This energy is invariant under M\oe bius
transformation and the "round" sphere is its only absolute minimum.Here
we sketch a proof
of the compactness property for a variant of it. The details will
appear elsewhere. \par\par
\bigskip

Given an embedding $f$ of $S^1$ into $R^3$ or $S^3$, it is possible
to associate to it an "energy"
defined as
$$
	E_0(f)=\int\int_{S^1\times S^1} \Big({1\over \vert x-y\vert ^2} -
 {1\over \vert f(x)-f(y)\vert^2}\Big) dxdy.
$$
where x and y are two points on the unit circle, $|x-y|$ is the length
of the cord and $|f(x)-f(y)|$
is the euclidean distance, in the case the ambient space is $R^3$, or
 the cord distance in the case
of  $S^3$. The term ${1\over \vert x-y\vert^2}$ is necessary to make
the integral
convergent and is usually called regularization in the literature.
 For a link e.g. $f\sqcup
g$  one adds to $E_0(f)$ and $E_0(g)$ an uregularized energy:
$$
	E_u(f;g)= \int\int_{S^1\times S^1}
 {1\over \vert f(x)-g(y)\vert^2} dxdy.
$$
where x and y lie on different circles.
 In [O'] O'Hara proved that if one adds to
$E_0(f)$ a term involving the geodesic curvature $\kappa (x)$ of
the knot and defines a new energy
$$
E_\lambda (f) = E_0 +\lambda \int_{S^1} (\kappa (x))^2dx
$$
with $\lambda >0$, the new functional has the compactness
property:\par \medskip
\proclaim The number of knot types that can be realized by an
 embedding $f$ with $ E_\lambda (f)<K$ is less than
$N(K)$ where $N(K)$ is some finite function of $K$.
\par \medskip
Freedman, He and Wang discovered that $E_0(f)$ is invariant
 under M\oe bius transformations of $S^3$ and
where able to prove, among other things, the two foundamental
properties:\par \par \medskip
1) The "round" circle is the only absolute minimum of $E_0(f)$ \par \par
2) $E_0(f)$ has the compactness property.\par \par \medskip
There has been some effort to find an analogue of the energy
for embeddings of $S^2$ and recently
Knuser and Sullivan [KS] and Auckly and Sadun [AS] have proposed
several candidates. The energy should
satisfy the following axioms, see [AS]:\par \par \medskip
1) $E_0(f)$ is invariant under M\oe bius transformations \par\par
2) $E_0(f)$ is bounded below.\par\par
3) $E_0(f)$ blows up when the embedding approaches a non-injective
 map.\par \par
4) Given two embeddings $f$ and $g$, for every $\epsilon$ there is
 an embedding $h$ in the isotopy class
of the connnected sum $f\#g $ such that
 $E_0(h)\leq E_0(f)+E_0(g)+\epsilon$.\par\par
5) There exists an unregularized energy for pairs of disjoint embeddings,
 $E_u(f;g)$, such that,
if we put $E(f\sqcup g)=E_0(f)+E_0(g)+E_u(f;g)$, $E(f\sqcup g)$ is
 invariant
 for M\oe bius
transformations.\par\par \medskip
The Energy of [KS] does not satisfy 5, that of [AS] is not known
to satisfy 2,
 in [D] we proposed an energy for
embeddings of $S^2$ into $R^n$ or $S^n$ that satisfies all the five
 axioms.\par
Moreover we showed that the round sphere is the only absolute
minimum for $E_0(f)$.
 Here we prove
that, if we add to $E_0(f)$ an additional term similar to the one of [O'],
our new functional has the compactness property for embeddings into
$R^4$ or $S^4$, which is
the only meaningful case.\par Our definition can be easily extended
to give an energy for
embeddings of surfaces of any genus. There is some evidence that,
 in the case of tori, the minimum is
given by the Clifford torus.
\bigskip
\centerline {\bf Definition of the Energy}
\medskip
The energy $E_0(f)$ was defined in [D] as:
$$
	\int\int_{S^2\times S^2} \Big({1\over \vert x-y\vert ^2} -
 {d_f(x)d_f(y)\over \vert f(x)-f(y)\vert^2}\Big)^2 dxdy.
$$
and there we sketched the proof that it satisfies the five
properties quoted before and also that
$E_0(f)=0$ if and only if $f$ is the embedding of the round $S^2$.
 Here we shall define a new energy:
$$
	E_\lambda(f)=
E_0(f)+ \lambda \Big (\int_{f(S^2)}(|H|^2-2K)^p dvol\Big
)^{1/p}{{Vol(f(S^2))}^{1/q}}
$$
where $H$ is the mean curvature vector of the surface in $R^n$,
 $\lambda>0$, $1/p+1/q=1$ and $p>1$.
Remark that the second term is invariant for dilatations but
 not for M\oe bius tranformations anymore.
However by the Schwarz inequality:
$$
\Big (\int_{f(S^2)}{(|H|^2-2K)^p}dvol\Big )^{1/p}{{Vol(f(S^2))}^{1/q}}\geq
\int_{f(S^2)}{(|H|}^2-2K) dvol
$$
and the last term is invariant under M\oe bius tranformations, so $E_\lambda
 (\mu\circ f)$ is bounded
below when $\mu$ varies over all the M\oe bius transformations which do not
send a
point of $f(S^2)$ to infinity.  Now we sketch the proof of the compactness
property:
\par\bigskip
\proclaim Theorem: Let $f$ be an embedding of $S^2$ into $R^4$ with
$E_\lambda(f) \leq C$, then $f$
can belong to a finite number $N(C)$ of isotopy classes of knots only.
\par\bigskip
Proof: We break the proof in two lemmas. In the first we will use the
 second term of the energy
to get some local regularity, roughly Lemma 1 shows that there is a
uniform length scale in which the
surface is very close in the $C^1$ norm to its tangent plane.The second
lemma uses a theorem of
Teichmueller on conformal
moduli of planar domains to show that an upper bound on $E_\lambda (f)$
gives a lower bound for the
Kuiper selfdistance of the image. The theorem is a consequence of this
last fact.\par\medskip
 \proclaim
Lemma 1 :Given any $\delta>0$ we can cover $f(S^2)$ with
$N(\delta,\lambda ,C)$ 4-balls $B_i$ so that for
any $i$ and any $ a\in B_i\cap f(S^2)$, if $Tf(S^2)|_a$ denotes the
tangent plane in a, the
ortogonal projection $\pi :B_i\mapsto Tf(S^2)|_a$ is a diffeomorfism
with biLipschitz constant $1+\delta /2$,
in particular it is a $1+\delta$ quasiconformal homeomorphism.
We can also assume that the center of $B_i$ is on $f(S^2)$.
\par\medskip
 Proof: Since $E_0(f)$ is always nonnegative we must have:\par
$$
\Big (\int_{f(S^2)}(|H|^2-2K)^p dvol\Big )^{1/p}\leq
{C\over \lambda{Vol(f(S^2))}^{1/q}}$$\par
now it is easy to see that the 2p-th powers of the components of the
second foundamental form
of the surface are bounded by $(|H|^2-2K)^p$ and so the left hand side
of the inequality is the
the square of the $L^{2p}$ norm of the second foundamental form.
The latter is also the derivative of
the Gauss map T from the surface into the Grassmannian of two planes
in $R^4$, it follows that the
inequality above gives a bound on the $L^{1,2p}$ norm squared of the
Gauss map.
We use the Sobolev embedding of $L^{1,2p}$ into $C ^{1/q}$, the
space of H\oe lder functions
of exponent $1/q$, to show that for any $\epsilon>0$ we have:
$$
|T(a)-T(b)|\leq
{A \Big ({C\over {\lambda (Vol(f(S^2)))^{1/q}}}\Big)^{1/2}|a-b|^{1/q}}
$$
here $a$ and $b$ are points on the surface, the vertical bars denote
the distances in the
Grassmannian and in the surface and $A$ is the constant of the Sobolev
embedding.
This estimate implies the lemma by an easy gemetric
argument.Q.E.D. \par\smallskip
In particular this tells us that we can cover $f(S^2)$
with a finite number (depending only on $C$)
of balls in such a way that in each of these balls the
 surface can be approximated by a disjoint union
 of planes, we will call them sheets.\par\medskip
\proclaim Lemma 2 : $E_\lambda (f)<C $ implies that,
 modulo a
refinement of the covering which will increase the
number of balls by a function of $C$ only, we can
 assume that each ball contains only one sheet.
\par \medskip
Proof: Let $a$ be the centre of a ball $B_i$ and $r$ be
 its radius.Let $\epsilon$ be the distance
between $a$ and a sheet in $B_i$ not containing it, if
 it exists. We shall prove that if $\epsilon$
is too small, the energy cannot be less than $C$.
Let $b$ be the point realizing the distance $\epsilon$
and let $U_a$ and $U_b$ be the sheets containing
$a$ and $b$ respectively; if $\epsilon$ is small enough
we can assume that the inner radius of both
sheets is at least $r/2$, say. Let $V_a$ and $V_b$ be
smaller disks in the induced metric on the
sheets whose radius is $mr$ where $m<1/2$ will be
determinated later.We want to prove that the
integral:
$$
	\int_{x \in {f^{-1}(V_a)}} \int_{y\in {f^{-1}(V_b)}}
\Big({1\over \vert x-y\vert ^2} -
 {d_f(x)d_f(y)\over \vert f(x)-f(y)\vert^2}\Big)^2 dxdy.
$$
goes to infinity when $\epsilon $ goes to zero.
First of all by the triangle inequality and substitution
of variables we have:
$$
\Big|{\int_{x \in {f^{-1}(V_a)}} \int_{y\in {f^{-1}(V_b)}}
\Big({1\over \vert x-y\vert ^2} -
 {d_f(x)d_f(y)\over \vert f(x)-f(y)\vert^2}\Big)^2 dxdy}
\Big|^{1/2} \geq
$$
$$
\geq \Big| \int_{V_a}\int_{V_b}
{1\over \vert \xi -\eta \vert ^4}d\xi d\eta \Big|^{1/2}-\Big|
\int_{f^{-1}(V_a)} \int_{f^{-1}(V_b)} {1\over \vert x-y\vert ^4}dxdy
\Big|^{1/2}
$$
where $\xi=f(x)$ and $\eta=f(y)$.\par
The proof of Lemma 1 shows that the first integral
in the right hand side is close to
the integral over the projection of $V_a$ and $V_b$ onto
their tangent planes at $a$ and $b$
respectively. The latter is greater than the integral over
two planar disks of radius $mr/(1+\delta)$
embedded in $R^4$ at Hausdorff distance less than $\epsilon$.
Here the radius of the disks does not depend
on $\epsilon$.The evaluation of this integral and
the proof that it goes to infinity when
$\epsilon \rightarrow 0$ is an exercice in calculus which we
willingly leave to the reader.\par
Now we prove that the second integral is bounded.
First remark that the $1+\delta$ quasi-conformality of the
projection onto the tangent spaces implies that
the annuli $U_a-V_a$ and $U_b-V_b$ have conformal modulus
very close to $log (2/m)$, and
that, since $f$ is  conformal, the same is true of their
preimages in $S^2$.\par
Modulo a M\oe bius transformation in the domain we can
assume that $f(0)=a$ and $f(\infty)=b$.
Now if $m$ is chosen small enough a theorem of Teichmueller
on extremal domains (see [LV], p.57)
says that
$$
{{inf_{x\in {S^2-f^{-1}(U_a)}}|a-x|} \over
{sup_{y\in {f^{-1}(V_a)}}|a-y|}}\geq D(m)
$$
where as usual the vertical bars denote the spherical distance in $S^2$.
 Note that D(m) does not
depend on $\epsilon$. The same estimate holds for $f^{-1}(U_b)$ and it
 is easy to see that these
two together imply that the integral
$\int_{f^{-1}(V_a)} \int_{f^{-1}(V_b)} {1\over \vert x-y\vert
^4}dxdy$ is bounded independently of $\epsilon$.
This ends the proof of Lemma 2 and gives
the compactness property.
Q.E.D.
\bigskip
\centerline{References}
\par \medskip
\frenchspacing
\item{[AS]} D.Auckly, L. Sadun, A family of M\oe bius
invariant 2-Knot energies, preprint
\item{[D]}  S. Demichelis, Energie des surfaces plong\'ees invariante par
transformation de M\oe bius, preprint
\item{[FHW]} M.H. Freedman,Z.X. He, Z. Wang, M\oe bius energy of
knots and unknots, Ann. of. Math.,139 (1994),1-50.
\item{[KS]} R.Knuser,J.Sullivan, M\oe bius Energies for Knots
and Links, preprint.
\item{[LV]} O. Lehto, K.I. Virtanen, Quasikonforme Abbildungen,
Springer-Verlag, 1965
\item{[O']} J. O'Hara, Energy of a knot, Topology, 30 (1991), 241-247
\bigskip
\bigskip
Stefano Demichelis \par
Dipartimento di matematica \par
Universita di Pavia  \par
27100 Pavia (Italie) \par
\end